\documentclass[11pt]{amsart}

\title[Quantum cohomology of $\cnums ^{2}/\Zmodthree$
and Hurwitz-Hodge integrals]{The orbifold quantum cohomology of $\cnums ^{2}/\Zmodthree$
and Hurwitz-Hodge integrals} 
\date{\today}
\address{
Dept of Math, Univ. of British Columbia,
Vancouver, BC, Canada 
}
\email{jbryan@math.ubc.ca}
\address{
Dept of Math, California Institute of Technology, Pasadena, CA
}
\email{graber@caltech.edu}
\address{
Dept of Math, Princeton Univ., Princeton, NJ 
}
\email{rahulp@math.princeton.edu}

\author{J. Bryan, T. Graber, and R.
Pandharipande} 

\date{\today}

\usepackage{diagrams}

\usepackage{amsmath}
\usepackage{amsmath,amsthm,amsfonts}

\newcommand{\C}{{\mathbb C}}
\newcommand{\Q}{{\mathbb Q}}
\newcommand{\cnums} {{\mathbb C}}          % complex numbers
\newcommand{\qnums} {{\mathbb Q}}		% rationals

\renewcommand{\P}{\mathbb{P}}
\newcommand{\M}{\overline{M}}
\renewcommand{\H}{\overline{H}}

\newcommand{\X}{\mathcal{X}}
\newcommand{\orb}{\mathit{orb}}
\newcommand{\util}{x}
\newcommand{\Ytil}{\widetilde{Y}}
\newcommand{\Xtil}{\widetilde{X}}
\newcommand{\F}{\widetilde{F}}
\renewcommand{\O}{\mathcal{O}}
\newcommand{\oh}{{\mathcal{O}}}
\newcommand{\Cbar}{\overline{C}}
\newcommand{\omegabar}{\overline{\omega}}
\newcommand{\pibar}{\overline{\pi}}
\newcommand{\E}{\mathbb{E}}
\newcommand{\rk}{\operatorname{rk}}

\newcommand{\Zmodthree}{{\mathbb{Z}_{3}}}
\newcommand{\Zmodtwo}{{\mathbb{Z}_{2}}}
\newcommand{\Zmodn}{{\mathbb{Z}_{n}}}

\newtheorem{thm}{Theorem}[section]

\newtheorem{cor}[thm]{Corollary}

\newtheorem{lemma}[thm]{Lemma}
\newtheorem{prop}[thm]{Proposition}
\newtheorem{proposition}[thm]{Proposition}

\newtheorem{conjecture*}{Conjecture}

\begin{document}
%\maketitle 
\begin{abstract}
Let $\Zmodthree$ act on $\cnums^2$ by non-trivial opposite characters.
Let $\X =[\cnums ^{2}/\Zmodthree]$ be the orbifold quotient, and let
$Y$ be the unique crepant resolution.  We show the equivariant genus 0
Gromov-Witten potentials $F^{\X}$ and $F^{Y}$ are equal after a change
of variables --- verifying the Crepant Resolution Conjecture for the
pair $(\X,Y)$. Our computations involve Hodge integrals on trigonal
Hurwitz spaces which are of independent interest. In a self contained
Appendix, we derive closed formulas for these Hurwitz-Hodge integrals.

\end{abstract}
\maketitle 
%\thanks{The authors are supported by NSERC, NSF, and \dots .}   

%\markboth{???}  {???}
%\renewcommand{\sectionmark}[1]{}

%\tableofcontents
%\pagebreak

\section{Introduction} 

The Crepant Resolution Conjecture predicts the Gromov-Witten
theory of a Gorenstein orbifold $\X $ is equivalent to the Gromov-Witten theory
of any crepant resolution $Y$. The conjecture
was originally formulated in physics by Zaslow and Vafa
\cite{Zaslow,Vafa-orbifold-numbers} and
subsequently 
in mathematics by Chen and Ruan \cite{Chen-Ruan}. 
A
precise statement of the general conjecture is given in
\cite{Bryan-Graber}.

One impediment to understanding the Crepant Resolution Conjecture is
the dearth of non-trivial examples where the full Gromov-Witten theory
(even in genus 0) of $\X$ and $Y$ has been computed. In
\cite{Bryan-Graber}, the genus 0 (equivariant) Crepant Resolution
Conjecture is verified in the cases
\begin{equation*} 
(\X ,Y) = ( \cnums
^{2}/\Zmodtwo,T^{*}\P ^{1}), \ \ (\X ,Y)=
(\operatorname{Sym}^{d}\cnums ^{2},\operatorname{Hilb}^{d}\cnums
^{2}).
\end{equation*}
These examples, while highly non-trivial, are limited in their
ability to exhibit many of the features of the general conjecture. In
particular, since the Picard numbers are 1, the
change of variables  has a restricted
form.

Our main result is the proof of the
equivariant genus 0 Crepant Resolution
Conjecture for the orbifold $\X = [\cnums ^{2}/\Zmodthree]$ with
unique crepant resolution $Y$. Here, the Picard number is 2, and
we see a more complicated transformation  taking place. Our
computations involve new integrals of Hodge classes over
trigonal Hurwitz spaces.

\subsection{Notation}
Let $\Zmodthree \subset SU (2)$ act on $\cnums ^{2}$ via the
standard representation of $SU (2)$. Let 
\[
\omega=e^{2\pi i/3}. 
\]
We identify $\Zmodthree$ with $\left\{1,\omega ,\omegabar \right\}$.
The $\Zmodthree$-action on $\cnums^2$ is
$$\omega\cdot (x,y) = (\omega x, \omegabar y).$$

Let $\mathcal{X}=[\cnums ^{2}/\Zmodthree]$ be the quotient stack with 
 coarse moduli space $X$. 
The singular variety $X$ admits a unique crepant resolution $$Y\to
X.$$ The exceptional divisor is a chain of two rational curves $E_{1}$
and $E_{2}$.  The action of the torus 
$$T=\cnums ^{*}\times \cnums^{*}$$ on 
$\cnums ^{2}$ commutes with the $\Zmodthree$-action  and induces 
$T$-actions on $\X $ and $Y$.

The potential  $F^{Y}$ is the
generating function for  equivariant genus 0 Gromov-Witten
invariants of $Y$:
\[
F^{Y}=\sum _{\beta=d_1[E_1]+d_2[E_2] }\ 
\sum _{n_{0},n_{1},n_{2}\geq 0}
\left\langle
1^{n_{0}}C_{1}^{n_{1}}C_{2}^{n_{2}} \right\rangle_{\beta
}^{Y} 
\frac{y_{0}^{n_{0}}}{n_{0}!}\frac{y_{1}^{n_{1}}}{n_{1}!}
\frac{y_{2}^{n_{2}}}{n_{2}!} 
q_{1}^{d_{1}}q_{2}^{d_{2}}.
\]
The first sum ranges over effective curve classes
$\beta$. 
The classes $C_i \in H_T^*(Y)$ are defined as
\begin{eqnarray*}
C_1 &= &  -\frac{2}{3} [E_1] - \frac{1}{3} [E_2] \\
C_2 & = &-\frac{1}{3} [E_1] - \frac{2}{3} [E_1].
\end{eqnarray*} 
The images of the $C_i$ in $H^*(Y)$ are Poincare dual to the proper
transforms of the images of the two coordinate axes in $\C^2$ but the 
equivariant lifts here are chosen to make them dual to the $[E_i]$ with 
respect to the equivariant intersection form.
The Gromov-Witten
invariants $\left\langle \cdot \right\rangle^Y_{\beta }$ are 
 multilinear functions on $H_T^*(Y)$ taking values in
$$H^{*}_{T} (pt)=\qnums [t_{1},t_{2}].$$  The unstable terms,
where $d_{1}=d_{2}=0$ and $n_{0}+n_{1}+n_{2}<3$, are defined to be zero.

%In \cite{Bryan-Graber}, orbifold Gromov-Witten invariants are extended
%to include a notion of degree for curves classes in twisted
%sectors. The corresponding Gromov-Witten potential has new quantum
%parameters labeled by $u$'s. 

For the orbifold $\X =[\cnums ^{2}/\Zmodthree]$, we
have a basis for $H^{*}_{\orb, T} (\X )$, the $T$-equivariant {\em orbifold}
cohomology of $\X $, given by classes $\{1,D_{1},D_{2} \}$ 
corresponding to the elements
$\{1,\omega ,\omegabar  \}$ of $\Zmodthree$. The potential
 $F^{\X }$ generates
the equivariant genus 0 orbifold Gromov-Witten invariants of $\X $:
\[
F^{\X }=
%\sum _{\beta }
%\sum _{\begin{smallmatrix} n_{0},n_{1},n_{2}\geq 0\\
%d_{1},d_{2}\geq 0  \end{smallmatrix}} 
\sum _{n_{0},n_{1},n_{2}\geq 0}
\left\langle
1^{n_{0}}D_{1}^{n_{1}}D_{2}^{n_{2}} \right\rangle^{\X } 
\frac{x_{0}^{n_{0}}}{n_{0}!}\frac{x_{1}^{n_{1}}}{n_{1}!}
\frac{x_{2}^{n_{2}}}{n_{2}!}. 
\]
The bracket $\langle \cdot \rangle ^{\X
}$ denotes the equivariant degree 0, genus 0  orbifold Gromov-Witten
invariant of  $\X $.  Again, we set unstable terms (those with fewer than three insertions) equal to zero.

\subsection{Results} The main result of the paper is the complete
computation of the potential functions $F^{\X }$ and $F^{Y}$.
The computations verify the Crepant Resolution
Conjecture for the pair $(\X,Y)$.

\begin{thm}\label{thm: formula for FY}
The equivariant genus 0 Gromov-Witten potential of $Y $ is:
\begin{align*}
F^{Y}=&\quad \frac{y_{0}^{3}}{18t_{1}t_{2}}-\frac{y_{0}}{3}
(y_{1}^{2}+y_{1}y_{2}+y_{2}^{2})\\
&+\frac{1}{3} (t_{1}+2t_{2})\frac{y_{1}y_{2}^{2}}{2}+\frac{2}{3} (2t_{1}+t_{2})\frac{y_{1}^{3}}{6}
+\frac{1}{3} (2t_{1}+t_{2})\frac{y_{1}^{2}y_{2}}{2}+\frac{2}{3} (t_{1}+2t_{2})\frac{y_{2}^{3}}{6}\\
&+ (t_{1}+t_{2})\sum _{d=1}^{\infty }\frac{1}{d^{3}}\left[
(e^{y_{1}}q_{1})^{d}+ (e^{y_{2}}q_{2})^{d}+
(e^{y_{1}+y_{2}}q_{1}q_{2})^{d} \right].
\end{align*}
\end{thm}

\begin{thm}\label{thm: formula for FX}
The equivariant genus 0 Gromov-Witten potential of $\X $ is:
\begin{align*}
F^{\X }=&\quad \frac{1}{18t_{1}t_{2}}x_{0}^{3}+ \frac{1}{3}x_{0} \util_{1} \util _{2} +\frac{1}{18}t_{1} \util _{1}^{3}+\frac{1}{18}t_{2} \util _{2}^{3}\\
&+ \frac{(t_{1}+t_{2})}{2}\sum _{g=2}^{\infty }\frac{(-1)^{g-1}
A_{g}}{(g+2)!}\frac{1}{3} \Big[ (\util _{1}+\util _{2})^{g+2}+ (\omega
\util _{1}+\omegabar \util _{2})^{g+2} \\
&  \hspace{200pt}+ (\omegabar \util _{1}+\omega \util_{2})^{g+2} \Big],
\end{align*}
where 
%\[
%%\util _{i}=u_{i}+x_{i}\quad \quad 
%\omega=e^{2\pi i/3}. 
%\]
%and 
the rational numbers $A_{g}$ are determined by:
\[
A (u)=\sum _{g=1}^{\infty }A_{g}\frac{u^{g-1}}{(g-1)!} =
\frac{1}{\sqrt{3}}\tan \left(\frac{u}{\sqrt{12}}+\frac{\pi }{6}
\right).
\]
\end{thm}

Geometrically, $A_{g}$ arises as the integral of
the Hodge class $\lambda _{g-1}$ over any connected component of the
Hurwitz scheme of curves in $\M _{g}$ which admit a cyclic triple
cover of $\P ^{1}$. In the Appendix, which is written to be
self-contained, we prove $A_{g}$ is independent of the choice of
component (Proposition~\ref{prop: integral is independent of
component}) and is given by the above formula
(Proposition~\ref{prop: generating function for Z/3 hodge
integrals}). We also prove formulas for related 
trigonal Hurwitz-Hodge integrals (Propositions~\ref{prop: generating function
for S3 hodge integrals}).

The series $F^{Y}$ converges at
$q_{i}=\omega$, in particular, the change of variables
\begin{align}
y_{0}&=x_{0}\label{eqn: change of vars first}\\
y_{1}&=\frac{i}{\sqrt{3}} (\omega \util _{1}+\omegabar \util _{2})\label{eqn: change of vars second}\\ 
y_{2}&=\frac{i}{\sqrt{3}} (\omegabar  \util _{1}+\omega\util _{2})\label{eqn: change of vars third}\\
q_{i}&=\omega \label{eqn: change of vars last}
\end{align}
is well-defined.
\begin{thm}\label{thm: FX=FY under change of variables}
After      the above change of variables, 
$$F^{\X }=F^{Y}$$ as power series in $x_{0}$, $\util _{1}$, and $\util
_{2}$ up to unstable terms. Hence, the equivariant genus 0 Crepant Resolution
Conjecture holds for $(\X,Y) $.
\end{thm}

\begin{cor}\label{cor: QH(X)=QH(Y)}
The equivariant 
quantum cohomology rings
$QH^{*}_{T,\orb } (\X )$ and $QH^{*}_{T} (Y)$ are isomorphic
after the above change of variables.
\end{cor}

The correct definition of quantum cohomology for Gorenstein orbifolds
requires the notion of quantum parameters in the twisted sector, see
\cite{Bryan-Graber}.

\subsection{DuVal singularities} Let $G\subset SU (2)$ be a finite
subgroup. Let 
\[
\X =[\cnums ^{2}/G]
\]
be the orbifold quotient, and let
$Y$ be the unique crepant resolution of the DuVal singularity $X$.

By the McKay correspondence, the cohomology of $Y$ has a natural basis
indexed by irreducible representations of $G$, where the trivial
representation $\cnums $ corresponds to the identity in $H^{0} (Y)$
and non-trivial representations $R$ correspond to classes in
$H^{2} (Y)$. The orbifold cohomology of $\X $ has a
natural basis indexed by conjugacy classes of $G$, where the trivial
conjugacy class $(e)$ corresponds to the identity in $H_{\orb }^{0}
(\X)$ and non-trivial conjugacy classes $(g)$ correspond to
classes in $H_{\orb }^{2} (\X)$.

We speculate that the potential functions $F^{\X}$ and $F^{Y}$ are
identified by the change of variables:
\begin{align*}
y_{\cnums }&=x_{(e)},\\
y_{R} &= \frac{1}{|G|}\sum _{g\in G} (\chi _{\rho } (g)-2)^{1/2}\, \chi
_{R} (g) \, x_{(g)},\\
q_{R}&=\omega ^{n_{R}}.
\end{align*}
Here $\rho $ is the standard representation of $G\subset SU (2)$ on
$\cnums ^{2}$, $\omega $ is a primitive $|G|$-th root of unity, and
$n_{R} $ is the coefficient of $R$ in the representation corresponding
the the longest root of the associated Dynkin diagram.

The above change of variables specializes to equations~\eqref{eqn:
change of vars first}--\eqref{eqn: change of vars last} for the case
of $\cnums ^{2}/\Zmodthree $.

\subsection{Acknowledgments} The authors are grateful to R. Cavalieri,
H. Esnault, and E.  Viehweg for helpful conversations.  

J.B. was supported by NSERC, T.G. was supported by the NSF and the
Sloan foundation, and R.P. was supported by the NSF and the Packard
foundation.  The research was partially pursued at the AMS summer
institute in algebraic geometry in Seattle, the Banff International
Research Station, and the Instituto Superior T\'ecnico in Lisbon.

\section{The Gromov-Witten invariants of $Y$} We
compute the equivariant genus 0 Gromov-Witten invariants of
$Y$ via an
equivariant embedding into a Calabi-Yau threefold $\Ytil$ for which the
Gromov-Witten invariants have been previously computed,

Consider the threefold $\Xtil \subset \cnums ^{4}$ given by the equation
\[
xy=z (z-s) (z+s).
\]
$\Xtil $ admits a small resolution $\Ytil \subset \cnums ^{4}\times \P
^{1}\times \P ^{1}$ given by the closure of the graph of the rational
map
\begin{align*}
\Xtil &\dasharrow \P ^{1}\times \P ^{1}\\
(x,y,z,s)&\mapsto (x:z),(x:z (z-s))
\end{align*}
see \cite{BKL,Katz}.

The surfaces $X$ and $Y$ are isomorphic to the 
subvarieties of $\Xtil $ and $\Ytil $ defined by 
$s=0$. This construction is $T$ equivariant under the action
\begin{align*}
(x,y,z,s)&\mapsto (t_{1}^{3}x,t_{2}^{3}y,t_{1}t_{2}z,t_{1}t_{2}s)\\
(u_{1}:v_{1}),(u_{2}:v_{2})&\mapsto
(t_{1}^{2}u_{1}:t_{2}v_{1}),(t_{1}u_{2}:t_{2}^{2}v_{2})
\end{align*}
 
The exceptional set of the resolution $\Ytil \to \Xtil $ consists of
two rational curves $E_{1}\cup E_{2}$ meeting in a point
$p_{1}$. $\Ytil $ is a Calabi-Yau threefold and the normal bundle of
$E_{i}\subset \Ytil $ is $\O (-1)\oplus \O (-1) $.

The Gromov-Witten invariants of $\Ytil $ were computed by
Bryan-Katz-Leung (see Proposition~2.10 of \cite{BKL}). For
\[
\beta =d_{1}E_{1}+d_{2}E_{2}\neq 0,
\]
the genus 0 invariants of $\Ytil $ are given by
\[
\left\langle \;  \right\rangle^{\Ytil }_{\beta }=\begin{cases}
\frac{1}{d^{3}}&\text{if $(d_{1},d_{2})= (d,d)$, $(d,0)$, or $(0,d)$}\\
0&\text{otherwise.}
\end{cases}
\]

The normal bundle $Y\subset \Ytil $ is trivial with the $T$-action
for which 
\[
c_{1} (N_{Y/\Ytil })=t_{1}+t_{2}.
\]
The 0 point invariants
of $\Ytil $ can be computed in terms of the 0 point invariants of
$Y$ as follows.
\begin{align*}
\left\langle \;  \right\rangle^{\Ytil }_{\beta }&=\int _{[\M _{0,0} (\Ytil ,\beta )]^{vir}}1\\
&=\int _{[\M _{0,0} (Y,\beta )]^{vir}}\frac{1}{e (R^{\bullet }\pi _{*}f^{*} (N_{Y/\Ytil }))}\\
&=\frac{1}{t_{1}+t_{2}}\int _{[\M _{0,0} (Y,\beta )]^{vir}}1\\
&=\frac{1}{t_{1}+t_{2}}\left\langle \;  \right\rangle^{Y}_{\beta }
\end{align*}
where $\pi :\mathcal{C}\to \M _{0,0} (Y,\beta )$ and $f:\mathcal{C}\to
Y$ are the universal curve and universal map respectively.  Combining
the above with the divisor and point axioms, we see that the $\beta
\neq 0$ part of $F^{Y}$ is given by:
\[
 (t_{1}+t_{2})\sum _{d=1}^{\infty }\frac{1}{d^{3}}\left[
(e^{y_{1}}q_{1})^{d}+ (e^{y_{2}}q_{2})^{d}+
(e^{y_{1}+y_{2}}q_{1}q_{2})^{d} \right].
\]

To finish the proof of Theorem~\ref{thm: formula for FY}, we must
to compute the $\beta =0$ terms of $F^{Y}$. These consist solely of
three point invariants given by triple intersections:
\[
\left\langle \gamma _{1},\gamma _{2},\gamma _{3} \right\rangle = \int
_{Y}\gamma _{1}\cup \gamma _{2}\cup \gamma _{3}
\]
which we define by localization and in general depend on the choice of
the equivariant lifts of $\gamma _{i}$.

The surface $Y$ has three $T$ fixed points $p_{0},p_{1},p_{2}$. The
$T$-invariant curve $E_{1}$ connects $p_{0}$ and $p_{1}$, and  
the $T$-invariant curve $E_{2}$
connects $p_{1}$ and $p_{2}$. The
weights of the $T$-action on $T_{p_{i}}Y$ can be easily computed from
our explicit description of $Y$ and are given by
\[
(3t_{1},-2t_{1}+t_{2}),\quad (2t_{1}-t_{2},-t_{1}+2t_{2}),\quad
(t_{1}-2t_{2},3t_{2})
\]
for $T_{p_{0}}Y$, $T_{p_{1}}Y$, and $T_{p_{2}}Y$ respectively.

The basis $\{C_{1},C_{2} \}$ is dual to the basis
$\{E_{1},E_{2} \}$. We can choose a lift of the $T$-action on $Y$ to
$L_{i}=\O (C_{i})$ such that the weights of the $T$-action on
$L_{1}|_{p_{0}}$, $L_{1}|_{p_{1}}$, $L_{1}|_{p_{2}}$ are 
\[
-2t_{1},\quad -t_{2},\quad -t_{2}
\]
respectively, and the weights of the $T$-action on
$L_{2}|_{p_{0}}$, $L_{2}|_{p_{1}}$, $L_{2}|_{p_{2}}$ are 
\[
-t_{1},\quad -t_{1},\quad -2t_{2}
\]
respectively. We can then compute by localization:
\[
\left\langle 1,1,1 \right\rangle = \frac{1}{3t_{1}t_{2}},\quad
\left\langle 1,1,C_{1} \right\rangle=0, \quad \left\langle 1,1,C_{2}
\right\rangle=0,
\]
\[
\left\langle 1,C_{1},C_{1} \right\rangle=-\frac{2}{3},\quad
\left\langle 1,C_{2},C_{2} \right\rangle=-\frac{2}{3},\quad
\left\langle 1,C_{1},C_{2} \right\rangle=-\frac{1}{3},
\]
\[
\left\langle C_{1},C_{1},C_{1} \right\rangle=\frac{2}{3}
(2t_{1}+t_{2}), \quad \left\langle C_{1},C_{1},C_{2} \right\rangle
=\frac{1}{3} (2t_{1}+t_{2}),
\]
\[
\left\langle C_{2},C_{2},C_{2} \right\rangle=\frac{2}{3}
(2t_{2}+t_{1}), \quad \left\langle C_{2},C_{2},C_{1} \right\rangle
=\frac{1}{3} (2t_{2}+t_{1}),
\]

completing the proof of Theorem~\ref{thm: formula for FY}.

\section{The Orbifold Gromov-Witten invariants of $X$} 
The cubic terms of $F^{\X}$ can be computed directly.  The
 higher degree terms are expressed here as trigonal Hurwitz-Hodge integrals
and computed in the Appendix.

The inertia stack $I\X$ has three components corresponding to the
three elements $\left\{1,\omega ,\omegabar \right\}$ of $\Zmodthree$. Each
component is contractable and so the graded vector space
\[
H^{*}_{\orb } (\X) = H^{*} (I\X)
\]
has a canonical basis $\left\{1,D_{1},D_{2} \right\}$ corresponding to
the three components. Moreover, the grading for the twisted sectors is
shifted by two:
\[
1\in H^{0}_{\orb } (\X) \quad \text{and} \quad D_{i}\in H^{2}_{\orb } (\X).
\]

The invariant $\left\langle 1^{n_{0}}D_{1}^{n_{1}}D_{2}^{n_{2}}
\right\rangle^\X$ is defined to be the integral
\[
\int _{\M _{0,n_{0}+n_{1}+n_{2}} (\X ,0)} \prod
_{i=1}^{n_{0}}{\text {ev}}_{i}^{*} (1)\prod_{ i=n_{0}+1}^{n_{0}+n_{1}}\text{ev}_{i}^{*}
(D_{1})\prod_{ i=n_{0}+n_{1}+1}^{n_{0}+n_{1}+n_{2}}{\text {ev}}_{i}^{*} (D_{2}).
\]

By the usual point axiom, 
%\cite{??} 
$\left\langle
1^{n_{0}}D_{1}^{n_{1}}D_{2}^{n_{2}} \right\rangle^\X=0$ if $n_{0}>0$ and
$n_{0}+n_{1}+n_{2}>3$. Moreover, if $n_{1}+n_{2}>0$, then there must
be stacky points of the domain curves of the twisted stable
maps. 
Consequently, the maps must factor through $B\Zmodthree\subset \X$. 

Consider $\left\langle D_{1}^{n_{1}}D_{2}^{n_{2}} \right\rangle^{\X}$
where $n_{1}+n_{2}>3$. Since the maps factor through $B\Zmodthree$, we can
rewrite the integral in terms of stable maps to $B\Zmodthree$:
\begin{multline*}
\left\langle D_{1}^{n_{1}}D_{2}^{n_{2}} \right\rangle^{\X }= \\
\int
_{[\M _{0,n_{1}+n_{2}} (B\Zmodthree)]^{vir}} e (R^{1}\pi _{*}f^{*} (L_{\omega
}\oplus L_{\omegabar })) \prod _{i=1}^{n_{1}}{\text {ev}}_{i}^{*} (D_{1})\prod
_{i=n_{1}+1}^{n_{1}+n_{2}}\text{ev}_{i}^{*} (D_{2}),
\end{multline*}
where $$\pi :\mathcal{C}\to \M _{0,n_{1}+n_{2}} (B\Zmodthree)$$ is the universal
curve and $$f:\mathcal{C}\to B\Zmodthree$$ is the universal map.
The normal bundle of $B\Zmodthree\subset
\X $ is the sum of the line bundles
$L_{\omega }\oplus L_{\omegabar }$ 
 determined by the
$\Zmodthree$-representations where $\omega \in \Zmodthree$ acts by multiplication by
$\omega $ and $\omegabar $ respectively.

Concretely, $\M _{0,n_{1}+n_{2}} (B\Zmodthree)$ may be thought of as
parameterizing curves $\Cbar $ equipped with a $\Zmodthree$-action for which
the quotient map is a cover $p:\Cbar \to C$ of a $n_{1}+n_{2}$ marked
genus 0 curve $C$ ramified over the marked points and possibly the
nodes of $C$. The integral $\left\langle D_{1}^{n_{1}}D_{2}^{n_{2}}
\right\rangle^{\X} $ is possibly non-zero only on the components of $\M
_{0,n_{1}+n_{2} } (B\Zmodthree)$ where $p:\Cbar \to C$ is ramified over all the
marked points with monodromy $\omega $ around the first $n_{1}$ points
and $\omegabar$ around the last $n_{2}$ points.

Consider the diagram of universal structures: 
\[
\begin{diagram}
\overline{\mathcal{C}}&\rTo&pt\\
\dTo^{p}&&\dTo\\
\mathcal{C} &\rTo ^{f}&B\Zmodthree\\
\dTo^{\pi }&&\\
\M _{0,n_{1}+n_{2}} (B\Zmodthree)&&
\end{diagram}
\]
Let $\pibar :\overline{\mathcal{C}}\to \M _{0,n_{1}+n_{2}} (B\Zmodthree)$ be
the composition $\pi \circ p$ and let
\[
\E ^{\vee }=R^{1}\pibar _{*}\O 
\]
be the dual Hodge bundle. By the Riemann-Hurwitz formula, $\E ^{\vee
}$ is a bundle of rank
\[
g=n_{1}+n_{2}-2.
\]
The action of $\omega \in \Zmodthree$ on $\overline{\mathcal{C}}$ induces an
action of $\omega $ on $\E ^{\vee }$. This gives a decomposition of $\E^\vee$
into
eigenbundles 
\[
\E ^{\vee }=\E ^{\vee }_{1}\oplus \E ^{\vee }_{\omega }\oplus \E
^{\vee }_{\omegabar }.
\]
(Note that our convention throughout is that 
$\E^{\vee }_\omega$ is the $\omega$ eigenbundle of $\E^\vee$ and not
the dual of $\E_\omega$.)

A chase through the definitions shows that 
\[
R^{1}\pi _{*}f^{*} (L_{\omega })=\E ^{\vee }_{\omegabar }, \quad R^{1}\pi
_{*}f^{*} (L_{\omegabar })=\E ^{\vee }_{\omega }.
\]
Moreover, $\E ^{\vee }_{1}=0$ is empty since 
$$\E^\vee_1 =R^{1}\pi _{*}\O $$ and $\pi
$ is a family of genus 0 curves.

Let $\M ^{\sigma }_{0,n_{1}+n_{2}}$ be the component of $\M
_{0,n_{1}+n_{2}} (B\Zmodthree)$ on which
\[
\prod _{i=1}^{n_{1}}\text{ev}_{i}^{*} (D_{1})\prod
_{i=n_{1}+1}^{n_{1}+n_{2}}\text{ev}_{i}^{*} (D_{2})
\]
is possibly non-zero. We can identify $\M ^{\sigma }_{0,n_{1}+n_{2}}
(B\Zmodthree)$ with the Hurwitz scheme $\H ^{\sigma }_{g} ((3)^{g+2})$ defined
in the Appendix.

So we have 
\[
\left\langle D_{1}^{n_{1}}D_{2}^{n_{2}} \right\rangle^{\X }=\int _{[\M
^{\sigma }_{0,n_{1}+n_{2}} (B\Zmodthree)]}e (\E ^{\vee }_{\omega }\oplus \E ^{\vee
}_{\omegabar })
\]
where $e$ is the $T$-equivariant Euler class. 

Since $\E ^{\vee }_{\omega }\oplus \E ^{\vee }_{\omegabar }$ has rank
$g=n_{1}+n_{2}-2$ and $\M _{0,n_{1}+n_{2}}^{\sigma } (B\Zmodthree)$ has
dimension $n_{1}+n_{2}-3$, $\left\langle D_{1}^{n_{1}}D_{2}^{n_{2}}
\right\rangle^{\X}$ is a linear function of the equivariant parameters
$t_{1}$ and $t_{2}$.

\begin{lemma}\label{lem: t1+t2 divides <D1n1D2n2>}
For $n_{1}+n_{2}>3$, $\left\langle D_{1}^{n_{1}}D_{2}^{n_{2}}
\right\rangle^{\X}$ is a multiple of $t_{1}+t_{2}$.
\end{lemma}
\textsc{Proof:} It suffices to prove that $\left\langle D_{1}^{n_{1}}D_{2}^{n_{2}} \right\rangle^{\X } =0$ for 
\[
t_{1}=-t_{2}=t.
\] 
Let 
\[
r_{1}=\rk \E ^{\vee }_{\omega }, \quad r_{2}=\rk \E ^{\vee }_{\omegabar }.
\]
Then,
$$e (\E ^{\vee }_{\omega }\oplus \E ^{\vee }_{\omegabar }) = (-1)^{r_{2}}\left(t^{g}+t^{g-1}c_{1} (\E ^{\vee }_{\omega }\oplus
\E _{\omega })+\dotsb +c_{g} (\E ^{\vee }_{\omega }\oplus \E_{\omega }) \right).$$
%\begin{align*}
%e (\E ^{\vee }_{\omega }\oplus \E ^{\vee }_{\omegabar })&=e (\E ^{\vee }_{\omega })e (\E ^{\vee }_{\omegabar })\\
%&= (t^{r_{1}}+t^{r_{1}-1}c_{1} (\E ^{\vee }_{\omega })+\dotsb
%+c_{r_{1}} (\E ^{\vee }_{\omega }))\cdot ((-t)^{r_{2}}+ (-t)^{r_{2}-1}c_{1}
%(\E ^{\vee }_{\omegabar })+\dotsb +c_{r_{2}} (\E ^{\vee }_{\omegabar
%}))\\
%&= (-1)^{r_{2}} (t^{r_{1}}+t^{r_{1}-1}c_{1} (\E ^{\vee }_{\omega })+\dotsb +c_{r_{1} }(\E ^{\vee }_{\omega }))\cdot (t^{r_{2}}+t^{r_{2}-1}c_{1} (\E _{\omegabar })+\dotsb +c_{r_{2}} (\E _{\omegabar }))\\
%&= (-1)^{r_{2}}\left(t^{g}+t^{g-1}c_{1} (\E ^{\vee }_{\omega }\oplus
%\E _{\omegabar })+\dotsb +c_{g} (\E ^{\vee }_{\omega }\oplus \E
%_{\omegabar }) \right)
%\end{align*}
The Lemma then follows from a $\Zmodthree$-version of Mumford's relation:
\begin{proposition}\label{prop: G-version of Mumford's relation}
Let $\pi :C\to B $ be a flat family of prestable curves with
the action of a finite group $G$. Let $\omega _{\pi }$ be the relative
dualizing sheaf and let $\E =\pi _{*}\omega _{\pi }$ be the Hodge
bundle. Let
\[
\E = \oplus _{\rho}\E _{\rho }
\]
be the decomposition of summands corresponding to the irreducible
representations of $G$. Then
\[
c (\E _{\rho }\oplus \E ^{\vee }_{\rho})=1\in H^*(B,\Q).
\]
\end{proposition}
\textsc{Proof:} The following argument is known to experts and is
referred to by Mumford in \cite{Mumford}, but since it does not seem to be
written down, we include it for the benefit of the reader.  We may
assume that $B$ is smooth, proper, and that the boundary divisor
$D\subset B$ over which $C$ is singular has normal crossings.  The
Lemma follows from the decomposition of $R^1\pi_*\C$ into eigensheaves
for the natural action of $G$.  Over $B-D$ we have the standard
sequence
$$ 0 \to \E^\vee \to R^1\pi_*\C\otimes \oh_B \to \E \to 0$$ which
admits an extension over all of $B$ where the middle term is
interpreted globally as 
\[
V=R^1\pi_*[\oh_C \rTo^{d} \omega_{C/B}].
\]
Moreover, the Gauss-Manin connection extends to a connection over all
of $B$ with logarithmic poles along $D$ whose polar part is nilpotent
(page 130, \cite{Griffiths}).  Because the Gauss-Manin connection
respects the decomposition of $R^1\pi_*\C$ into eigenbundles, it
follows that this extension does so as well.  Thus, after splitting
$V$ into eigenbundles, we get a sequence
$$ 0 \to \E_{\rho}^\vee \to V_\rho \to \E_\rho \to 0 $$ and on
$V_\rho$ we have a log connection with nilpotent residue.  In
\cite[Appendix B]{Esnault-Viehweg}, it is shown how to use a
connection with log poles to compute the Atiyah class (and hence the
Chern classes) of a bundle.  Because the formula for the Chern classes
is in terms of the eigenvalues of the residue of the connection, and
these all vanish in our situation, it follows that the Chern classes
of $V_\rho$ all vanish.

\qed

\vspace{+10pt}

By Lemma \ref{lem: t1+t2 divides <D1n1D2n2>}, after setting
\[
t_{1}=t_{2}=t,
\]
we obtain
\begin{align*}
\left. \left\langle D_{1}^{n_{1}}D_{2}^{n_{2}} \right\rangle^{\X } \right|_{t_1=t_2=t} &=\int
_{\M ^{\sigma }_{0,n_{1}+n_{2}} (B\Zmodthree)}e (\E ^{\vee }_{\omega }\oplus \E
^{\vee }_{\omegabar })\\
&=t\int _{\M ^{\sigma }_{0,n_{1}+n_{2}} (B\Zmodthree)}c_{g-1} (\E ^{\vee })\\
&=t (-1)^{g-1}\int _{\H ^{\sigma }_{g} ((3)^{g+2})}\lambda _{g-1}\\
&=t (-1)^{g-1}A_{g},
\end{align*}
where the last equality is well defined by Proposition~\ref{prop:
integral is independent of component} and the values of $A_{g}$ are
given by Proposition~\ref{prop: generating function for Z/3 hodge
integrals}.

For
$
g=n_{1}+n_{2}-2>1,
$
we conclude
\[
\left\langle D_{1}^{n_{1}}D_{2}^{n_{2}} \right\rangle^{\X }=\begin{cases}
\frac{t_{1}+t_{2}}{2} (-1)^{g-1}A_{g} &\text{for $n_{1}\equiv  n_{2}\mod 3$}\\
0&\text{for $n_{1}\not \equiv  n_{2}\mod 3$}.
\end{cases}
\]
Let
$$F^\X= F^\X_{\text{cubic}} + \widehat{F}^\X,$$
where $F^\X_{\text{cubic}}$ consists of all the cubic terms.
Then,
\begin{align*}
\widehat{F}^{\X} &= \frac{t_{1}+t_{2}}{2}\sum _{g=2}^{\infty } \sum _{\begin{smallmatrix} n_{1}+n_{2}=g-2\\
n_{1}\equiv n_{2}\mod 3 \end{smallmatrix}}  (-1)^{g-1}A_{g} \frac{x_{1}^{n_{1}}}{n_{1}!}  \frac{x_{2}^{n_{2}}}{n_{2}!}\\
&=\frac{t_{1}+t_{2}}{2}\sum _{g=2}^{\infty } (-1)^{g-1}\frac{A_{g}}{(g-2)!} \sum _{\begin{smallmatrix} n_{1}+n_{2}=g-2\\
n_{1}\equiv n_{2}\mod 3 \end{smallmatrix}} \binom{g-2}{n_{1}}x_{1}^{n_{1}}x_{2}^{n_{2}}\\
&=
\quad \frac{t_{1}+t_{2}}{2}\sum _{g=2}^{\infty }
(-1)^{g-1}\frac{A_{g}}{(g-2)!}\frac{1}{3}\Big[ (x_{1}+x_{2})^{g-2}+
(\omega x_{1}+\omegabar x_{2})^{g-2}\\
& \ \ \ \ \ \ \ \ \ \ \ \ \  \ \ \ \ \
\ \ \ \ \ \ \ \ \ \ \ \ \ \ \ \ \ \ \ \ \ \ \ \ \ \ \ \ \ \ \ \ \ \ \ \ \ \ \ \ \ \ \ \ 
+ (\omegabar x_{1}+\omega
x_{2})^{g-2} \big].
\end{align*}
To finish the proof of Theorem~\ref{thm: formula for FX}, we must
remains to compute $F^\X_{\text{cubic}}$.

%\smargin{I wrote this last bit here in a really hurried fashion. We need to work on it a bit}
By the monodromy condition, the only non-vanishing 3-point invariants are
\[
\left\langle 1^{3} \right\rangle^{\X },\quad \left\langle 1D_{1}D_{2}
\right\rangle^{\X },\quad \left\langle D_{1}^{3} \right\rangle^{\X
},\quad \left\langle D_{2}^{3} \right\rangle^{\X }.
\]

The moduli space for the first invariant is just $\X$ itself, so it is
a trivial localization calculation.  Each of the other invariants 
is an integral over a moduli space consisting of a single point with a
$\Zmodthree$ automorphism group. In the first case, the corresponding cover is
connected and genus 0, and in the last two cases, $\Cbar $ is the
elliptic curve with an order 3 automorphism (the two non-trivial
automorphisms determining the two different cases). Thus the
invariants are all $1/3$ times the appropriate weight, namely
\[
\frac{e (H^{1} (\Cbar ,\O )_{\omega }\oplus H^{1} (\Cbar ,\O
)_{\omegabar })}{e (H^{0} (\Cbar ,\O )_{\omega }\oplus H^{0} (\Cbar
,\O )_{\omegabar })}
\]
where the subscript indicates the eigenspace for the action of $\omega
$. These are easily computed by first principles or by the holomorphic
Lefschetz formula. We get
\[
\left\langle 1^{3} \right\rangle^{\X }=\frac{1}{3t_{1}t_{2}},\quad
\left\langle 1D_{1}D_{2} \right\rangle^{\X }=\frac{1}{3},\quad \left\langle
D_{1}^{3} \right\rangle^{\X }=\frac{t_{1}}{3},\quad \left\langle D_{2}^{3}
\right\rangle^{\X }=\frac{t_{2}}{3}.
\]
The proof of Theorem~\ref{thm: formula for FX} is complete. \qed

\vspace{+10pt}

Trigonometric evaluations of Hodge integrals over hyperelliptic Hurwitz
spaces \cite{Faber-Pandharipande-logarithmic} play a basic role in the
Crepant Resolution Conjecture for $\cnums^2/\Zmodtwo$ studied in
\cite{Bryan-Graber}.  Hodge integrals over trigonal Hurwitz spaces
arise in the study of $\cnums^2/\Zmodthree$. For $n\geq 4$, the Crepant
Resolution Conjecture for $\cnums^2/\Zmodn$ (where $\Zmodn$ acts
on the first factor via the standard representation $\rho$ and the
second factor via the dual representation $\rho^\vee$) predicts simple
evaluations of certain Chern classes of $\E_\rho \oplus
\E_{\rho^\vee}$ on Hurwitz spaces of $\Zmodn$-covers.

\section{Checking the series agree}

We make the substitutions given by equations \eqref{eqn: change of
vars first}--\eqref{eqn: change of vars last} into $F^{Y} $ and compare
with $F^{\X}$. The terms of homogeneous degree $-2$ and degree $0$ in
$t_{1}$ and $t_{2}$ are easily checked to agree with the corresponding
terms in $F^{\X}$.

The remaining terms are linear in $t_{1}$ and $t_{2}$. 
Hence, agreement for the specializations
$t_{1}+t_{2}=0$ and $t_{1}-t_{2}=0$ implies full agreement. The case of $t_{1}+t_{2}=0$ is
straightforward. 

Let $\F^{Y}$ and $\F^{\X}$ denote the $t$ linear term of
$F^{Y}|_{t_{1}=t_{2}=t}$ and $F^{\X}|_{t_{1}=t_{2}=t}$
respectively. We need to prove that after making the substitution
\eqref{eqn: change of vars first}--\eqref{eqn: change of vars last},
$\F ^{Y}$ and $\F ^{\X}$ agree as power series in $\util _{1}$ and
$\util _{2}$ up to terms of degree less than or equal to
two. Equivalently, we must check that the third partial
derivatives of $\F^{Y}$ and $\F ^{\X}$ agree.
\begin{align*}
\F^{\X}=& \frac{1}{18} (\util _{1}^{3}+\util _{2}^{3})\\  & +\quad \sum
_{g=2}^{\infty }\frac{(-1)^{g-1} A_{g}}{(g+2)!}\frac{1}{3}\left[
(\util _{1}+\util _{2})^{g+2}+ (\omega \util _{1}+\omegabar \util _{2})^{g+2}+ (\omegabar\util _{1}+\omega
\util_{2})^{g+2} \right]\\
=&\sum _{g=1}^{\infty }\frac{(-1)^{g-1}
A_{g}}{(g+2)!}\frac{1}{3}\left[ (\util _{1}+\util _{2})^{g+2}+ (\omega
\util _{1}+\omegabar \util _{2})^{g+2}+ (\omegabar \util _{1}+\omega \util_{2})^{g+2} \right].
\end{align*}
Differentiation yields formulas for the partial derivatives (denoted by subscripts) :
\begin{align*}
 \F_{111} ^{\X}
=& \frac{1}{3}\left[A (-\util _{1}-\util _{2})+ A (-\omega \util
_{1}-\omegabar \util _{2})+A (-\omegabar \util
_{1}-\omega \util _{2}) \right],\\
\F_{112} ^{\X } =& \frac{1}{3}\left[A (-\util _{1}-\util _{2})+\omega
A (-\omega \util _{1}-\omegabar \util _{2})+\omegabar A (-\omegabar \util _{1}-\omega \util _{2}) \right].
\end{align*}
Similarly, for $Y$, we have
\begin{align*}
\F ^{Y}&= \frac{1}{2}y_{1}y_{2}^{2}+\frac{1}{3}y_{1}^{3}+
\frac{1}{2}y_{2}y_{1}^{2}+\frac{1}{3}y_{2}^{3}\\
 &\quad +  2\sum _{d=1}^{\infty
}\frac{1}{d^{3}}\left[ (e^{y_{1}}q_{1})^{d}+ (e^{y_{2}}q_{2})^{d}+
(e^{y_{1}+y_{2}}q_{1}q_{2})^{d} \right]\\
&=\left(\tfrac{i}{\sqrt{3}} \right)^{3}\left[\frac{\util _{1}^{3}}{6}+ \frac{\util _{2}^{3}}{6}- \util _{1}^{2}\util _{2}- \util _{2}^{2}\util _{1} \right]\\
&\quad +2\sum _{d=1}^{\infty }\frac{1}{d^{3}}\left[ \left(\omega
e^{\frac{i}{\sqrt{3}}\left(\omega \util _{1}+\omegabar \util
_{2} \right)} \right)^{d}+ \left(\omega
e^{\frac{i}{\sqrt{3}}\left(\omega \util _{2}+\omegabar \util
_{1} \right)} \right)^{d} +\left(\omegabar e^{-\frac{i}{\sqrt{3}}\left(\util _{1}+\util _{2} \right)}
\right)^{d} \right]
\end{align*}
In the variables
\begin{align*}
\xi &=\util _{1}+\util _{2},\\
\xi _{\omega }&=\omega \util _{1}+\omegabar \util _{2},\\
\xi _{\omegabar }&=\omegabar \util _{1}+\omega \util _{2},
\end{align*}
the partial derivative $\partial ^{3}/\partial x_{1}^{3}$ is:
\[
\F ^{Y }_{111}=\left(\tfrac{i}{\sqrt{3}} \right)^{3}\left(
1+\frac{2\omega e^{\frac{i}{\sqrt{3}}\xi _{\omega }}}{1-\omega
e^{\frac{i}{\sqrt{3}}\xi _{\omega }}} +\frac{2\omega
e^{\frac{i}{\sqrt{3}}\xi _{\omegabar }}}{1-\omega
e^{\frac{i}{\sqrt{3}}\xi _{\omegabar }}} -\frac{2\omegabar e^{-\frac{i}{\sqrt{3}}\xi }}{1-\omegabar e^{-\frac{i}{\sqrt{3}}\xi }} \right)
\]
Applying the identity
\[
\frac{2e^{2i\theta }}{1-e^{2i\theta }} = -i\tan\left(\theta
+\tfrac{\pi }{2} \right)-1,
\]
we obtain
\begin{align*}
\F ^{Y }_{111}&=\left(\tfrac{1}{3\sqrt{3}} \right)\left\{-\tan
\left(\tfrac{1}{\sqrt{12}} \xi _{\omega }+\tfrac{5\pi }{6} \right) -\tan
\left(\tfrac{1}{\sqrt{12}} \xi _{\omegabar }+\tfrac{5\pi }{6}
\right) +\tan \left(-\tfrac{1}{\sqrt{12}} \xi +\tfrac{7\pi }{6} \right)
  \right\}\\
&= \frac{1}{3}\left( A (-\xi _{\omega })+A (-\xi _{\omegabar })+A (-\xi ) \right)\\
&=\F ^{\X}_{111}.
\end{align*}
A similar computation verifies $\F ^{\X}_{112}=\F ^{Y}_{112}$. The
identities $\F ^{\X}_{122}=\F ^{Y}_{122}$ and $\F ^{\X}_{222}=\F
^{Y}_{222}$ are obtained by  symmetry in the indices.
Theorem~\ref{thm: FX=FY under change of variables} is proved.

The six unstable terms of $F ^{\X}$ can be assigned values
(expressed in terms of trilogarithms, dilogarithms, and logarithms) by imposing the
equality
 $F^{\X}=F^{Y}$. It would be interesting to give a geometric
interpretation of these unstable values.

\appendix \section{Degree 3 Hurwitz Hodge Integrals}\label{Appendix}

\subsection{} Consider the moduli spaces
$\H _g(\mu^1,\ldots,\mu^n)$ of connected, genus $g$, degree
$3$ admissible covers of an unparameterized $\mathbf{P}^1$.  We label
the monodromy conditions $\mu^i$ in degree 3 by the size of the
largest part of the associated partition.  There are two natural maps:
$$\epsilon:\H _g(\mu^1,\ldots,\mu^n)\rightarrow \M _g,$$
$$\pi: \H _g(\mu^1,\ldots,\mu^n) \rightarrow \M _{0,n},$$
well defined if $g\geq 2$ and $n\geq 3$ respectively.

\subsection{} We will primarily be interested in the moduli spaces
$\H _g((3)^{g+2})$ for $g\geq 1$.  Consider a covering
\[
[f:\Cbar  \rightarrow (C,p_1,\ldots,p_{g+2})]\in \H _g((3)^{g+2})
\]
where $(C,p_1,\ldots,p_n)$ is a stable, $n$-pointed, genus 0 curve.
Since all the monodromy conditions are 3-cycles, the covering $f$ {\em
must} be Galois with group $\Zmodthree$.

The monodromy around each ramification point determines a non-zero
element of the Galois group of $f$.  Hence, a canonical assignment
\[
\sigma_f: \{ p_1, \ldots, p_{g+2}\} \rightarrow
\text{Gal}(f)\setminus 0
\]
is determined by $f$.  Since
$\text{Gal}(f)\setminus 0$ has two elements, $\sigma_f$ defines a two
set partition of the markings,
\[
\{ p_1, \ldots, p_{g+2}\}= S_f \cup S'_f.
\]
The parity condition
\begin{equation}\label{red}
|S_f|=|S'_f| \mod 3
\end{equation}
must be satisfied by global monodromy considerations.

The connected components of $\H _g((3)^{g+2})$ are in
bijective correspondence with {\em unordered} partitions $S\cup S'$ of
the marking set satisfying the parity condition \eqref{red}.  Let
\[
\H ^\sigma_g((3)^{g+2})
\]
be the connected component corresponding to a partition $\sigma$ of
the marking set satisfying the parity condition \eqref{red}.

The total number $\gamma_g$ of connected components of
 $\H _g((3)^{g+2})$ is given by the following formula:
\[
\gamma_g= \frac{1}{2} \sum_{l=1-g \mod 3} \binom{g+2}{l}.
\]
The prefactor $1/2$ occurs since the set partition $\sigma$
is unordered.

\subsection{}
We calculate the  evaluations of $\lambda_{g-1}$ against the
components of the moduli space $\H _g((3)^{g+2})$.
For $g\geq 1$, let
\[
A^\sigma_g = \int_{\H ^\sigma_g((3)^{g+2})} \lambda_{g-1}.
\]

\begin{prop}\label{prop: integral is independent of component}
The integral $A^\sigma_g$ is independent of $\sigma$.
\end{prop}

Let $A_g$ be the common value of the evaluations of $\lambda_{g-1}$ over
the components of $\H ^\sigma_g((3)^{g+2})$.
Consider the generating function
$$A(u) = \sum_{g\geq 1} A_g \frac{u^{g-1}}{(g-1)!}.$$

\begin{prop} \label{prop: generating function for Z/3 hodge integrals} The generating function $A$ is determined by:
$$A(u)= \frac{1}{\sqrt{3}} \tan(\frac{u}{\sqrt{12}} + \frac{\pi}{6}).$$
\end{prop}

Let $A^\bullet_g$ denote the evaluation of $\lambda_{g-1}$ against
 the full moduli space $\H _g((3)^{g+2})$,
$$A^\bullet_g = \int_{\H ((3)^{g+2})} \lambda_{g-1}.$$
The relation
$$A^\bullet_g =\gamma_g \cdot A_g.$$
is a consequence of Propositions \ref{prop: integral is independent of component} and \ref{prop: generating function for Z/3 hodge integrals}.

\subsection{}

Our proofs of Propositions \ref{prop: integral is independent of component} and \ref{prop: generating function for Z/3 hodge integrals} require the study of
a closely related Hodge integral series.
For $g\geq 0$, let
$$B_g = \int_{\H _g((3)^{g+1}(2)^2)} \lambda_{g},$$
and let
$$B(u) = \sum_{g\geq 0} B_g \frac{u^g}{g!}.$$

\begin{prop} \label{prop: generating function for S3 hodge integrals} The generating function $B$ is determined by:
$$B(u)= \frac{1}{\sqrt{3}} \tan(\frac{u}{\sqrt{12}} + \frac{\pi}{3}).$$
\end{prop}

\subsection{}
We start by considering the initial values $B_0$ and $B_1$.
% With the exception of $B_1$,
% the higher values will be determined by geometric recursion.
% $$A_1 = \int_{\H _1((3)^{3})}1 = 1/3$$
% is a straightforward genus 1 Hurwitz number.
% Similarly
The first, 
$$B_0= \int_{\H _0((3)(2)^2)}1 = 1,$$
is a genus 0 Hurwitz number.
The second,
$$B_1= \int_{\H _1((3)^2(2)^2)}\lambda_1 = 2/3,$$
calculated by the following geometric argument.
Consider the map
$$\epsilon:\H _1((3)^2(2)^2)\rightarrow \M _{1,1}$$
obtained by marking the first triple ramification point.  By
definition,
$$B_1= \int_{\H _1((3)^2(2)^2)}\epsilon^*(\lambda_1).$$
However, on $\M _{1,1}$, $\lambda_1= \psi_1$.
Hence,
$$B_1=\int_{\H _1((3)^2(2)^2)}\epsilon^*(\psi_1) =  \int_{\H _1((3)^2(2)^2)} \psi_1.$$
The last equality is {\em not} formal, but rather proven geometrically
since the component of the admissible cover carrying
the first marking is {\em never} contracted by $\epsilon$.
Now, $\psi_1$ on ${\H _1((3)^2(2)^2)}$ is easily seen to be given by pull-back via $\pi$,
 $$\psi_1= \frac{1}{3}\pi^*(\psi_1).$$
By applying the boundary relation to $\psi_1$ on $\M _{0,4}$ and the degeneration
formula, we obtain the
answer.

\subsection{}

We now proceed to determine all the higher $B_g$. The method is a use of the WDVV relation in the
context of Hodge integrals over the moduli spaces of admissible covers. 
We will prove the following recursion for
$g\geq 2$,
\begin{equation}\label{ee}
B_{g-1}+ \sum_{g=h_1+h_2} 3\binom{g-2}{h_1} B_{h_1}B_{h_2} = 
\sum_{g=h_1+h_2} 6 \binom{g-2}{h_1-1} B_{h_1}B_{h_2}
\end{equation}
Since the left side contains the summand $3B_0B_g$ and the right side does not contain $B_g$,
all higher $B_g$ are determined.

For $g\geq 2$, consider the space of Hurwitz covers of
the {\rm rigid} line 
$$\H ^{r}_{g}((2)(2)(3)(3) (2)^2 (3)^{g-2}).$$
The dimension is $g+4$.
The ramification conditions are written as above to distinguish the first 4. 
Let $\xi$
be the class of a point on the rigid line.
We may consider the integral
$$C_g = \int_{\H ^{r}_{g}((2)(2)(3)(3) (2)^2 (3)^{g-2})} \lambda_g \cup \prod_{i=1}^4 
{\text {ev}}_i^*(p)$$
which fixes the positions of the first 4 ramification conditions.
We then may specialize the 4 points to be in WDVV configurations 
$$((2)(2)|(3)(3))  \ \ \text{and} \ \ ((2)(3)| (2)(3))$$
by breaking the rigid target.

We consider first the evaluation of $C_g$ via the configuration $((2)(2)|(3)(3))$. 
We must now distribute the remaining ramifications $(2)^2 (3)^{g-2}$ to either side.
We focus our attention on  the $(2)^2$.
\begin{enumerate}
\item[(i)] Both $(2)^2$ go to the left. By parity, the central partition over the node must be (1) or (3).
           If (1), then the resulting configuration must have a loop and is annihilated by $\lambda_g$.
           If (3), then the left moduli space is $\H _{h_1}((2)^4 (3)^{h_1})$ and the
right moduli space is $\H _{h_2}((3)(3)(3)^{h_2})$. However, the integrand distributes by
   $\lambda_{h_1}$ and $\lambda_{h_2}$ respectively. The dimension mismatch yields vanishing.

\item[(ii)] One $(2)$ goes to left and one $(2)$ goes to the right. By parity, the central partition
must be (2). Since the final configuration can not have a loop, only one possibility is allowed:
all $(3)^{g-2}$ are distributed to the right. The outcome is the term
$$2^2 \int_{\H _{1}((2)^4)} \lambda_1  \cdot B_{g-1}.$$
Here, one prefactor  of 2 comes from the initial choice of $(2)$ and one comes from
the degeneration formula.
The integral 
$$\int_{\H _{1}((2)^4)} \lambda_1 = 1/4$$
is easily evaluated.
\item[(iii)] Both $(2)^2$ go to the right. By parity, the central partition must by (1) or (3).
If (1), the the resulting configuration must have a loop and is annihilated by $\lambda_g$.
If (3), we obtain the sum
$$ \sum_{g=h_1+h_2} 3\binom{g-2}{h_1} B_{h_1}B_{h_2}.$$
\end{enumerate}

Next, we consider the evaluation of $C_g$ via the configuration $((2)(3)|(2)(3))$.
Since loops must be avoided, the only possibility for the central partition is $(3)$.
Hence, among the distributed ramifications, one $(2)$ must go to either side.
The outcome is 
$$\sum_{g=h_1+h_2} 2\cdot 3 \binom{g-2}{h_1-1} B_{h_1}B_{h_2}$$
completing the derivation of equation \eqref{ee}.
Here, a prefactor 2 comes from the initial choice of $(2)$ and a prefactor of 3 comes from
the degeneration formula.
The determination of the integrals $B_g$ is complete.

\subsection{} 

Multiplying equation \eqref{ee} by $u^{g-2}/ (g-2)!$ and summing over
all $g\geq 2$, we easily derive the following differential equation
for $B (u)$:
\[
B' +3B B''  = 6 (B')^{2}.
\]
With the initial conditions $B (0)=1$ and $B' (0)=2/3$, the above ODE
is uniquely solved by
\[
B (u)=\frac{1}{\sqrt{3}}\tan \left(\frac{u}{\sqrt{12}}+\frac{\pi }{3} \right)
\]
which proves Proposition \ref{prop: generating function for S3 hodge
integrals}.

\subsection{}

We now turn to the integrals $A^\bullet_g$. A similar Hurwitz Hodge WDVV argument yield the following
relation for $g\geq 1$,
\begin{equation}\label{ww}
\delta_{g,1}+ \sum_{g=h_1+h_2} 3 \binom{g-1}{h_1-1} A^\bullet_{h_1} B_{h_2} =
\sum_{g-1=h_1+h2} 2 \binom{g-1}{h_1} B_{h_1}B_{h_2}.
\end{equation}
Certainly equation \eqref{ww} determines all the integrals $A^\bullet_g$ 
from the integrals $B_g$.

For $g\geq 1$, consider the space of Hurwitz covers of
the {\rm rigid} line 
$$\H ^{r}_{g}((2)(2)(3)(3) (3)^{g-1}).$$
The dimension is $g+3$
The ramification conditions are written as above to distinguish the first 4. 
We may consider the integral
$$D_g = \int_{\H ^{r}_{g}((2)(2)(3)(3) (3)^{g-1})} \lambda_{g-1} \cup \prod_{i=1}^4 
{\text {ev}}_i^*(\xi)$$
which fixes the positions of the first 4 ramification conditions.
We then may specialize the 4 points to be in WDVV configurations 
$$((3)(3)|(2)(2))  \ \ \text{and} \ \ ((2)(3)| (2)(3))$$
by breaking the rigid target.

We consider first the evaluation of $C_g$ via the configuration $((3)(3)|(2)(2))$.
We must now distribute the remaining ramifications $(3)^{g-1}$ to either side.
By parity the central partition must be (1) or (3).
\begin{enumerate}
\item[(i)] If the central partition is (1) and at least one $(3)$ is distributed right, then 
 the resulting configuration must have two loops and is then annihilated by $\lambda_{g-1}$.
Hence, if the the central partition is (1), all $(3)^{g-1}$ must be distributed left.
The left moduli space is then $\H _{g-1}((3)^{g+1} (1))$
For $g\geq 3$, the map
$$\epsilon:\H _{g-1}((3)^{g+1} (1))\rightarrow \M _g$$
has 1-dimensional fibers and evaluated to 0 against any Hodge classes.
For $g=2$, the vanishing still holds by the 1-dimensional fibers of 
$$\epsilon:\H _{1}((3)^{3} (1))\rightarrow \M _{1,1}.$$
The only contribution comes when $g=1$. Then
configuration yields 
$$\delta_{g,1}.$$
\item[(ii)] If the central partition is (3), the outcome is the term
$$\sum_{g=h_1+h_2} 3 \binom{g-1}{h_1-1} A^\bullet_{h_1} B_{h_2}.$$
\end{enumerate}

Next, we consider the evaluation of $D_g$ via the configuration
$((2)(3)|(2)(3))$.  By parity, the only possibility for the central
partition is $(2)$.  The outcome is
$$\sum_{g-1=h_1+h2} 2 \binom{g-1}{h_1} B_{h_1}B_{h_2}$$ completing the
derivation of equation \eqref{ww}.  Here, a prefactor 2 comes from the
degeneration formula.  The determination of the integrals $A^{\bullet
}_g$ is complete.

%Some values (by hand, so not guaranteed):
%
%$$B_0=1, \ \ B_1=2/3,  \ \ B_2=2/3,$$
%$$A_1=1/3, \ \ A_2=2/3, \ \ A_3=10/27.$$

\subsection{}\label{subsec: derive A series, assuming components are equal}

Let 
\[
A^{\bullet } (u) = \sum _{g=1}^{\infty }A^{\bullet }_{g}\frac{u^{g-1}}{(g-1)!}.
\]

Multiplying equation \eqref{ww} by $u^{g-1}/ (g-1)!$ and summing over
$g\geq 1$, we easily derive the following relation for $A^{\bullet } (u)$ in
terms of $B (u)$:
\[
1+3A^{\bullet }B=2B^{2}
\]
or equivalently
\[
A^{\bullet }=\frac{2}{3}B - \frac{1}{3}B^{-1}.
\]

We will now prove Proposition~\ref{prop: generating function for Z/3
hodge integrals} \emph{assuming} Proposition~\ref{prop: integral is
independent of component}. We begin by finding a closed formula for
$\gamma _{g}$, the number of components of $\H _{g}
((3)^{g+2})$.

\begin{lemma}\label{lemma: formula for gamma}
The number of unordered set partitions $S\cup
S'=\{p_{1},\dotsc ,p_{g+2} \}$ satisfying $|S|\equiv |S'| \mod 3$ is
given by
\[
\gamma _{g}=\frac{1}{3} (2^{g+1}+ (-1)^{g}).
\]
\end{lemma}
\textsc{Proof:} Consider all unordered set partitions $S\cup
S'=\{p_{1},\dotsc ,p_{g+2} \}$ and let $\overline{S}\cup
\overline{S}'$ be the induced partition of $\{p_{1},\dotsc ,p_{g+1}
\}$. The partitions fall into three mutually exclusive possibilities:
\begin{enumerate}
\item[(i)] $|S|\equiv |S'| \mod 3$,
\item[(ii)] $|\overline{S}|\equiv |\overline{S}'| \mod 3$, or
\item[(iii)] neither equality holds.
\end{enumerate}

Set (i) has cardinality $\gamma _{g}$, set (ii) has cardinality
$2\gamma _{g-1}$, and set (iii) has a bijection with set (i) obtained by
moving $p_{g+2}$ from one set in the partition to the other.

Consequently, we obtain the following recursion for $\gamma _{g}$:
\[
2^{g+1}=2\gamma _{g}+2\gamma _{g-1}.
\]

The formula in the Lemma uniquely solves this recursion with the
initial condition $\gamma _{0}=1$. \qed

Now we \emph{assume} Proposition~\ref{prop: integral is
independent of component} holds so $A^{\bullet }=A_{g}\gamma
_{g}$. Then applying Lemma~\ref{lemma: formula for gamma}, we see 
\begin{align*}
A^{\bullet } (u)&=\sum _{g=1}^{\infty }A_{g}\gamma _{g}\frac{u^{g-1}}{(g-1)!}\\
&=\sum _{g=1}^{\infty }\frac{4}{3}A_{g}\frac{(2u)^{g-1}}{(g-1)!}-\frac{1}{3}A_{g}\frac{(-u)^{g-1}}{(g-1)!}\\
&=\frac{4}{3}A (2u)-\frac{1}{3}A (-u).
\end{align*}

To prove Proposition~\ref{prop: generating function for Z/3 hodge
integrals} (assuming Proposition~\ref{prop: integral is independent of
component}), we must verify that the series 
\begin{align*}
B (u)&=\frac{1}{\sqrt{3}}\tan \left(\frac{u}{\sqrt{12}}+\frac{\pi }{3} \right)\\
A (u)&=\frac{1}{\sqrt{3}}\tan \left(\frac{u}{\sqrt{12}}+\frac{\pi }{6} \right)\\
\end{align*}
satisfy the functional equation
\[
\frac{2}{3}B (u)-\frac{1}{3}B (u)^{-1}=\frac{4}{3}A (2u)-\frac{1}{3}A (-u).
\]
Let 
\[
x=\frac{u}{\sqrt{12}}+\frac{\pi }{3}
\]
Multiplying the functional equation by $3\sqrt{3}$, we get
\[
2\tan (x) - 3 \cot (x)=4\tan \left(2x-\frac{\pi }{2} \right)
+\tan\left(x-\frac{\pi }{2} \right).
\]
Applying the trigonometric identities
\[
\tan \left(\theta -\frac{\pi }{2} \right)=-\cot (\theta ), \quad \cot
(2\theta )=\frac{1}{2}\left(\cot (\theta )-\tan (\theta ) \right),
\]
the equality is easily seen to hold. We have proven:
\begin{lemma}\label{lem: prop(components equal)=>prop(generating fnc for Ag)}
Proposition~\ref{prop: integral is independent of
component} implies Proposition~\ref{prop: generating function for Z/3 hodge integrals}.
\end{lemma}

\subsection{} To prove Proposition~\ref{prop: integral is independent
of component}, we derive a set of recursions for the integrals
$A^{\sigma }_{g}$. These recursions, \emph{combined with}
the determination of $A^{\bullet }_{g}$, uniquely determine the values
of all the integrals $A^{\sigma }_{g}$. Since the
recursions are indeed satisfied when $A^{\sigma }_{g}=A_{g}$, the 
integrals $A^{\sigma }_{g}$
are independent of the component type and their values
are given by the generating function in Proposition~\ref{prop:
generating function for Z/3 hodge integrals}.

The method requires a WDVV equation for Hodge integrals on the
components of $\H _{g}((3)^{g+2})$.

Let $\sigma$ be a two set partition of the markings $\{p_1,\ldots,p_{g+1}\}$
satisfying the parity condition.
The integral $A^\sigma_g$ depends only on the length 2 partition
$$|S_\sigma| + |S'_\sigma|=g+2$$
as the geometry of the moduli space is symmetric under permutation of the markings.
Let
$$A^{l,l'}_g = A^\sigma_g$$
where $l+l'$ is the associated length 2 partition of $g+2$.

We must calculate all the integrals $A^{l,l'}_g$ where
$l+l'=g+2$ and
\[
l\equiv l' \mod 3.
\]
The constraints imply
\[
l\equiv 1-g \mod 3.
\]
In particular,
$$A_g^\bullet = \frac{1}{2} \sum_{l\equiv 1-g \; (3)} \binom{g+2}{l} \cdot A^{l,g+2-l}_{g}$$
where the prefactor $1/2$ corrects for the double counting since
\begin{equation}
\label{dfg}
A^{l,l'}_g= A^{l',l}_g
\end{equation}
correspond to the same class of components.

To simplify the notation,
will we often write $A^l_g$ for $A^{l,g+2-l}_g$. The equality
$$A_g^l = A_g^{g+2-l}$$
is obtained from \eqref{dfg}

\subsection{}
For $g\leq 3$, only a single length 2 partition of $g+2$ occurs
in each genus:
$$A_1^{0},\ \ A_2^{2}, \ \  A_3^{1}.$$
Hence, Proposition \ref{prop: integral is independent of component} is empty. 

%The integrals $A_g$ are determined for $g\leq 3$
%from the computation of $A^\bullet_g$ by
%$$A_1= A_1^\bullet,$$
%$$A_2= \frac{A_2^\bullet}{3},$$
%$$A_3=\frac{A_3^\bullet}{5}.$$ In particular, Proposition \ref{prop:
%generating function for Z/3 hodge integrals} is proven up to order
%$g=3$ by the computation in section~\ref{subsec: derive A series,
%assuming components are equal}.

\subsection{} Let $g\geq 4$ and assume Proposition \ref{prop: integral
is independent of component} is proven for all lower genera.  We will
now prove Proposition \ref{prop: integral is independent of component}
for genus $g$.

Let $\omega ,\omegabar $ denote the non-zero elements of
$\mathbb{Z}/3{\mathbb{Z}}$.  Let $2\leq l\leq g+1$ satisfy
\begin{equation} \label{qwq}
l=g+3-l \mod 3.
\end{equation}
Consider the connected component
$$\H _{g+1}(\omega ^l\omegabar ^{g+3-l})\subset \H
_{g+1}((3)^{g+3})$$ corresponding to the monodromy $\omega $ for the
first $l$ markings and $\omegabar $ for the last $g+3-l$.  Equation
\eqref{qwq} is the parity condition.

Let $p_1,p_{2}, q_1,q_2$ be the first and last two markings, and let
\[
\pi:\H _{g+1}(\omega ^l\omegabar ^{g+3-l}) \rightarrow \M
_{0,4}
\]
be the associated map.  Let
\[
E_{g+1}^{l} = \int_{\H _{g+1}(\omega ^l\omegabar ^{g+3-l})}
\lambda_{g-1} \cup \pi^{-1}(\xi)
\]
where $\xi$ is a class of a point in $\M _{0,4}$.

We may calculate $E_{g+1}^{l}$ by specializing $\xi$ to either of the 
two WDVV configurations
$$(p_1p_2|q_1q_2), \ \ \text{and} \ \  (p_1q_1|p_2q_2)$$
in $\M _{0,4}$.
The resulting equation is easily derived:
\begin{multline*}
\sum_{x-y\not\equiv  1 \; (3)} 3 \binom{l-2}{x} \binom{g+1-l}{y} A_{1+x+y}^{2+x+\phi(x,y)}
A_{1+(l-x)+(g+3-l-y)}^{l-x+\overline{\phi}(x,y)} = \\
\sum_{x-y\not\equiv  0 \; (3)} 3 \binom{l-2}{x} \binom{g+1-l}{y} A_{1+x+y}^{2+x+\theta(x,y)}
 A_{1+(l-x)+(g+3-l-y)}^{l-x+\overline{\theta}(x,y)}.
\end{multline*}
The functions $\phi$, $\overline{\phi}$, $\theta$, $\overline{\theta}$ 
are defined as follows:
\begin{align*}
 x-y\equiv 0 \mod 3: & \ \ \ \ \  \phi(x,y)=1, \ \overline{\phi}(x,y)=0 \\
x-y\equiv 1 \mod 3: &  \ \ \ \ \  \theta(x,y)=0, \ \overline{\theta}(x,y)=1 \\
x-y\equiv 2 \mod 3: &  \ \ \ \ \ \phi(x,y)=0, \ \overline{\phi}(x,y)=1 \\
&                \ \ \ \ \  \theta(x,y)=1, \ \overline{\theta}(x,y)=0.
\end{align*}
Let $\mathcal{E}_{g+1}^l$ denote the equation obtained from $E_{g+1}^l$.

No terms of $\mathcal{E}_{g+1}^l$ contain $A$-integrals of genus greater than $g$.
The {\em principal} terms of $\mathcal{E}_{g+1}^l$ are those
which contain $A$-integrals of genus $g$.
In fact, the principal terms of $\mathcal{E}_{g+1}^l$ occur only on the
left side and are simply
$$A_{g}^{l-2} + A_{g}^{l+1}.$$

We now study the full linear system of principal terms.
Let $\nu$ be the smallest
non-negative integer congruent to $(1-g)$ mod 3.
The set of $A$-integrals of genus $g$ is 
$$\{ A_g^\nu, \ A_g^{3+\nu}, \ A_g^{6+\nu}, \  \ldots,\ A_g^{g+2-\nu}\}.$$
To simplify notation, denote these $A$-integrals by the variables
$$x_i= A_g^{3i+\nu}.$$
We must solve for the variables
$$\{ x_0, \ldots, x_n\}$$
for $n=\frac{g+2-2\nu}{3}$.
Elementary considerations show the number of variables, $n+1$, is congruent to $(g+1)$ mod 2.

The set of principal terms of all the ${\mathcal E}_{g+1}^l$ equations is simply
\begin{equation}\label{xdfd}
\{ x_0+x_1,\ x_1+x_2,\ x_2+x_3, \ \ldots,\ x_{n-1}+x_n\}.
\end{equation}
These principal terms do not determine the variables: exactly one additional
independent equation is required.

If $g$ is odd, the equations $x_i=x_{n-i}$ from \eqref{dfg} provide an independent relation
since the number of variables then is even. 
If $g$ is even, the symmetry $x_i=x_{n-i}$ is redundant.

An additional linear equation is obtained from the completed
calculation of $A_g^\bullet$:
\begin{equation}\label{vvv}
\sum_{i=0}^n \binom{g+2}{3i+\nu} x_i =2 A_g^\bullet
\end{equation}
Equation \eqref{vvv} is independent of the principal terms
\eqref{xdfd} if and only if
$$
\delta _{g}=\sum_{i=0}^n \binom{g+2}{3i+\nu} (-1)^{3i+\nu}
$$
does not vanish.

\begin{lemma}The numbers $\delta _{g}$ are given by
\[
\delta _{g}=\begin{cases}
-2 (-3)^{g/2}&\text{if $g$ is even,}\\
0&\text{if $g$ is odd.}
\end{cases}
\]
\end{lemma}

\noindent \textsc{Proof:} Let $\omega= \exp(2\pi i/3)$, and note that
$\omega -1 =\sqrt{3}\exp (5\pi i/6)$. Define 
\[
\theta (x)=\begin{cases}
1&\text{if $x\equiv 0\mod 3$,}\\
0&\text{otherwise.}
\end{cases}
\]

\begin{align*}
\delta _{g}&=\sum _{i=0}^{n}\binom{g+2}{3i+\nu } (-1)^{3i+\nu }\\
&=\sum _{k=0}^{g+2}\binom{g+2}{k} (-1)^{k}\theta (1-g-k)\\
&=\sum _{k=0}^{g+2}\binom{g+2}{k} (-1)^{k}\frac{1}{3}\left(1+\omega ^{1-g-k}+\omegabar ^{1-g-k} \right)\\
&=\frac{1}{3} (1-1)^{g+2}+\frac{1}{3}\omega ^{-2g-1} (\omega -1)^{g+2}+\frac{1}{3}\omegabar ^{-2g-1} (\omegabar -1)^{g+2}\\
&=\frac{1}{3} (\sqrt{3})^{g+2}\left( e^{\frac{2\pi i (-2g-1)}{3}}e^{\frac{(2\pi i)5 (g+2)}{12}}+e^{ \frac{2\pi i (2g+1)}{3}}e^{\frac{(-2\pi i)5 (g+2)}{12}} \right)\\
&=3^{g/2}\left(e^{2\pi i\frac{2-g}{4}}+e^{2\pi i\frac{g-2}{4}} \right)\\
&=3^{g/2}\cdot \begin{cases}
0&\text{if $g$ is odd,}\\
2 (-1)^{g/2-1}&\text{if $g$ is even.}
\end{cases}
\end{align*}
\qed

We conclude the full set of component $A$-integrals is completely determined
by the following three conditions:
\begin{enumerate}
\item[(i)] the initial values for $g\leq 3$,
\item[(ii)] the equations ${\mathcal E}_{g+1}^l$ for $g\geq 4$.
\item[(iii)] the additional equation \eqref{vvv} obtained from $A_g^\bullet$.
\end{enumerate}
To complete the proof of Proposition \ref{prop: integral is
independent of component}, we must simply check the compatibility of
the proposed values for the component $A$-integrals
 with the conditions (i-iii).

Compatibility with (i) and (iii) has already been
checked. Compatibility with (ii) is equivalent to the following set of
relations for $A_g$: for every pair $(r,s)$ of non-negative integers
congruent mod 3 and not equal (0,0),
\begin{multline*}
\sum_{x-y\not\equiv  1 \; (3)}  \binom{r}{x} \binom{s}{y} A_{1+x+y}
A_{1+(r-x)+(s-y)} = \\
\sum_{x-y\not\equiv  0 \; (3)}  \binom{r}{x} \binom{s}{y} A_{1+x+y}
 A_{1+(r-x)+(s-y)}.
\end{multline*}

We define
\begin{align*}
\theta _{0,r,s}&=\sum_{x-y\equiv 0 \; (3)}  \binom{r}{x} \binom{s}{y} A_{1+x+y}
A_{1+(r-x)+(s-y)}  \\
\theta _{1,r,s}&=
\sum_{x-y\equiv  1 \; (3)}  \binom{r}{x} \binom{s}{y} A_{1+x+y}
 A_{1+(r-x)+(s-y)}.
\end{align*}

The above compatibility condition is equivalent to the condition that
$\theta _{0,r,s}=\theta _{1,r,s}$ for all $r\equiv s \mod 3$,
$(r,s)\neq (0,0)$. This is easily seen by subtracting the full sum
over $x$ and $y$ from both sides of the compatibility equation. Let
\[
\theta _{i} (v,w)=\sum _{
\begin{smallmatrix} 
r,s\geq 0\\
r\equiv s\; (3)
\end{smallmatrix}
}\theta
_{i,r,s}\frac{v^{r}}{r!}\frac{w^{s}}{s!}.
\]
We need to prove that 
\[
\theta _{0} (v,w)-\theta _{1} (v,w)=\frac{1}{9}.
\]
We expand $\theta _{i}$ and rearrange the sums:
\begin{align*}
\theta _{i} (v,w)&=\sum _{\begin{smallmatrix} r,s\geq 0\\r\equiv s \;
(3) \end{smallmatrix}} \sum _{\begin{smallmatrix} x,y\geq 0\\x\equiv
y+i\; (3) \end{smallmatrix}}\frac{r!}{x!(r-x)!}\frac{s!}{y!(s-y)!}A_{1+x+y}A_{1+r+s-x-y}\frac{v^{r}}{r!}\frac{w^{s}}{s!}\\
&=\sum _{\begin{smallmatrix} x,y\geq 0\\x\equiv y+i\; (3) \end{smallmatrix}} \sum _{\begin{smallmatrix} n,m\geq 0\\ n\equiv m-i\; (3) \end{smallmatrix}} \frac{(x+y)!}{x!\; n!}\frac{(n+m)!}{y!\; m!}A_{1+x+y}A_{1+n+m}\frac{v^{n+x}}{(x+y)!}\frac{w^{m+y}}{(n+m)!}\\
&=\sum _{\begin{smallmatrix} x,y\geq 0\\x\equiv y+i\;
(3) \end{smallmatrix}}\binom{x+y}{x}A_{1+x+y}\frac{v^{x}w^{y}}{(x+y)!}\sum
_{\begin{smallmatrix} n,m\geq 0\\n\equiv m-i\; (3) \end{smallmatrix}}\binom{n+m}{m}A_{1+n+m}\frac{v^{n}w^{m}}{(n+m)!}\\
&=Q_{i} (v,w)Q_{-i} (v,w)
\end{align*}
where
\[
Q_{i} (v,w)=\frac{1}{3} \left(\tilde{Q}_i (v,w)+\omegabar ^{i}\tilde{Q}_i (\omega v,\omegabar w)+\omega ^{i}\tilde{Q}_i (\omegabar v,\omega w) \right)
\]
and
\begin{align*}
\tilde{Q}_{i} (v,w)&=\sum _{x,y\geq 0}\binom{x+y}{x}v^{x}w^{y}\frac{A_{1+x+y}}{(x+y)!}\\
&=\sum _{k=0}^{\infty }A_{1+k}\frac{(v+w)^{k}}{k!}\\
&=A(v+w).
\end{align*}
Therefore
\begin{align*}
9 (\theta _{0}-\theta _{1}) &=9 ( Q_{0}^{2}-Q_{1}Q_{-1})\\
&= (A[0]+A[1]+A[2])^{2}- (A[0]+\omegabar A[1]+\omega A[2]) (A[0]+\omega A[1]+\omegabar A[2])\\
&=3 (A[0]A[1]+A[1]A[2]+A[2]A[0])
\end{align*}
where $A[k]$ is defined by:
\[
A[k] (v,w)=A (\omega ^{k}v+\omegabar ^{k}w).
\]

Now
\begin{align*}
A (u)&=\left(\frac{1}{\sqrt{3}} \right)\tan\left(\frac{u}{\sqrt{12}}+\frac{\pi }{6} \right)\\
&=\left(\frac{-i}{\sqrt{3}} \right)\frac{e^{i\left(u/\sqrt{12}+\pi /6 \right)}-e^{-i\left(u/\sqrt{12}+\pi /6 \right)}}{e^{i\left(u/\sqrt{12}+\pi /6 \right)}+e^{-i\left(u/\sqrt{12}+\pi /6 \right)}}\\
&=\left(\frac{i}{\sqrt{3}} \right)\frac{1+\omegabar e^{iu/\sqrt{3}}}{1-\omegabar e^{iu/\sqrt{3}}}.
\end{align*}
Let
\[
X=e^{iv/\sqrt{3}}\quad Y=e^{iw/\sqrt{3}},
\]
then
\[
A[k]=\left(\frac{i}{\sqrt{3}} \right)\frac{1+\omegabar X^{\omega ^{k}}Y^{\omegabar ^{k}}}{1-\omegabar X^{\omega ^{k}}Y^{\omegabar ^{k}}}
\]
Define
\[
\Theta ^{\pm }_{k}=1\pm \omegabar X^{\omega
^{k}}Y^{\omegabar ^{k}}
\]
then
\begin{align*}
-9 (\theta _{0}-\theta _{1})&=-3 (A[0]A[1]+A[1]A[2]+A[2]A[0])\\
&=\frac{\Theta ^{+}_{0}\Theta ^{+}_{1}}{\Theta ^{-}_{0}\Theta ^{-}_{1}}+\frac{\Theta ^{+}_{1}\Theta ^{+}_{2}}{\Theta ^{-}_{1}\Theta ^{-}_{2}}+\frac{\Theta ^{+}_{2}\Theta ^{+}_{0}}{\Theta ^{-}_{2}\Theta ^{-}_{0}} \\
&=\frac{\Theta ^{+}_{0}\Theta ^{+}_{1}\Theta ^{-}_{2}
+\Theta ^{+}_{0}\Theta^{-}_{1}\Theta ^{+}_{2}
+\Theta ^{-}_{0}\Theta ^{+}_{1}\Theta ^{+}_{2}}{\Theta ^{-}_{0}\Theta ^{-}_{1}\Theta ^{-}_{2}}
\end{align*}
Applying the relations $1+\omega +\omegabar =0$ and
$\omegabar ^{2}=\omega $, we compute:
\begin{align*}
\Theta _{0}^{-}\Theta _{1}^{-}\Theta _{2}^{-}&= (1-\omegabar XY) (1-\omegabar X^{\omega }Y^{\omegabar }) (1-\omegabar X^{\omegabar }Y^{\omega })\\
&=\omega (X^{-1}Y^{-1}+X^{-\omega }Y^{-\omegabar }+X^{-\omegabar }Y^{-\omega }).
\end{align*}
Similarly, we have
\begin{align*}
\Theta ^{+}_{0}\Theta ^{+}_{1}\Theta ^{-}_{2}&=
\omega (-X^{-1}Y^{-1}-X^{-\omega }Y^{-\omegabar }+X^{-\omegabar }Y^{-\omega })\\
\Theta ^{+}_{0}\Theta ^{-}_{1}\Theta ^{+}_{2}&=
\omega (-X^{-1}Y^{-1}+X^{-\omega }Y^{-\omegabar }-X^{-\omegabar }Y^{-\omega })\\
\Theta ^{-}_{0}\Theta ^{+}_{1}\Theta ^{+}_{2}&=
\omega (+X^{-1}Y^{-1}-X^{-\omega }Y^{-\omegabar }-X^{-\omegabar  }Y^{-\omega }),
\end{align*}
and so we conclude that
\[
-9 (\theta _{0}-\theta _{1})=-1
\]
as desired.

We've shown that the values of the integrals given by
Propositions~\ref{prop: integral is independent of component} and
\ref{prop: generating function for Z/3 hodge integrals} are indeed the
unique solution to the full set of recursions and so the proofs of
Proposition \ref{prop: integral is independent of component} and
\ref{prop: generating function for Z/3 hodge integrals} are
complete. \qed

\bibliography{mainbiblio}

\begin{thebibliography}{10}

\bibitem{Bryan-Graber}
Jim Bryan and Tom Graber.
\newblock The crepant resolution conjecture.
\newblock In preparation.

\bibitem{BKL}
Jim Bryan, Sheldon Katz, and Naichung~Conan Leung.
\newblock Multiple covers and the integrality conjecture for rational curves in
  {C}alabi-{Y}au threefolds.
\newblock {\em J. Algebraic Geom.}, 10(3):549--568, 2001.
\newblock Preprint version: math.AG/9911056.

\bibitem{Chen-Ruan}
Weimin Chen and Yongbin Ruan.
\newblock Orbifold {G}romov-{W}itten theory.
\newblock In {\em Orbifolds in mathematics and physics (Madison, WI, 2001)},
  volume 310 of {\em Contemp. Math.}, pages 25--85. Amer. Math. Soc.,
  Providence, RI, 2002.

\bibitem{Esnault-Viehweg}
H{\'e}l{\`e}ne Esnault and Eckart Viehweg.
\newblock Logarithmic de {R}ham complexes and vanishing theorems.
\newblock {\em Invent. Math.}, 86(1):161--194, 1986.

\bibitem{Faber-Pandharipande-logarithmic}
C.~Faber and R.~Pandharipande.
\newblock Logarithmic series and {H}odge integrals in the tautological ring.
\newblock {\em Michigan Math. J.}, 48:215--252, 2000.
\newblock With an appendix by Don Zagier, Dedicated to William Fulton on the
  occasion of his 60th birthday.

\bibitem{Griffiths}
Phillip Griffiths, editor.
\newblock {\em Topics in transcendental algebraic geometry}, volume 106 of {\em
  Annals of Mathematics Studies}, Princeton, NJ, 1984. Princeton University
  Press.

\bibitem{Katz}
Sheldon Katz.
\newblock Small resolutions of {G}orenstein threefold singularities.
\newblock In {\em Algebraic geometry: Sundance 1988}, pages 61--70. Amer. Math.
  Soc., Providence, RI, 1991.

\bibitem{Mumford}
David Mumford.
\newblock Towards an enumerative geometry of the moduli space of curves.
\newblock In {\em Arithmetic and geometry, Vol. II}, pages 271--328.
  Birkh\"auser Boston, Boston, Mass., 1983.

\bibitem{Vafa-orbifold-numbers}
Cumrun Vafa.
\newblock String vacua and orbifoldized {L}{G} models.
\newblock {\em Modern Phys. Lett. A}, 4(12):1169--1185, 1989.

\bibitem{Zaslow}
Eric Zaslow.
\newblock Topological orbifold models and quantum cohomology rings.
\newblock {\em Comm. Math. Phys.}, 156(2):301--331, 1993.

\end{thebibliography}
\bibliographystyle{plain}

\end{document}